\documentclass[11pt,dvips]{amsart}  
\usepackage{amssymb,amsmath,amsfonts}  
\usepackage[dvips]{graphicx,color}
\usepackage{cite}
\allowdisplaybreaks[4]
\definecolor{pgray}{gray}{0.8}

\newtheorem{theorem}{Theorem}[section]

\newtheorem{remark}[theorem]{Remark}
\numberwithin{equation}{section}

\title{\mbox{}}

\begin{document}
\begin{center}
{\bf \LARGE{
	Logarithmic evolutions on the incompressible Navier--Stokes flow
}}\\
\vspace{5mm}
Masakazu Yamamoto\footnote{e-mail : \texttt{mk-yamamoto@gunma-u.ac.jp}}\\Graduate School of Science and Technology, Gunma University
\end{center}
\maketitle
\vspace{-15mm}
%
\begin{abstract}
Through the asymptotic expansion, the large-time behavior of the incompressible Navier--Stokes flow in $n$-dimensional whole space is drawn.
In particular, the logarithmic evolution included in the flow velocity is the focus of attention.
When the components of velocity are ordered from major to minor according to the parabolic scales, the logarithmically evolving components appear in a certain pattern.
This work asserts that this pattern varies depending on the even-oddness of the space dimension.
This means that the parity of the nonlinear drift is different in even and odd dimensions.
The logarithmic term represents the nonlinear component of the phenomenon.
No symmetry of the initial condition is required to prove this fact.
In the preceding works, the expansion with the $2n$th order was already derived.
The assertion is derived by reexamining these works in detail.
\end{abstract}

\section{Introduction}
We study large-time behavior of the following initial value problem of incompressible Navier--Stokes equations in whole space:
\begin{equation}\label{NS}
 \left\{
\begin{aligned}
	&\partial_t u + (u\cdot\nabla) u = \Delta u - \nabla p,
	&
	t>0,~ x \in \mathbb{R}^n,\\
	&\nabla\cdot u = 0,
	&
	t>0,~ x \in \mathbb{R}^n,\\
	&u(0,x) = a (x),
	&
	x \in \mathbb{R}^n,
\end{aligned}
 \right.
\end{equation}
	%
%
where $n \ge 2$, and $u = (u^1,u^2,\ldots,u^n) (t,x) \in \mathbb{R}^n$ and $p = p(t,x) \in \mathbb{R}$ are unknown velocity and pressure, respectively, and $a = (a^1,a^2,\ldots,a^n) (x) \in \mathbb{R}^n$ is given initial velocity satisfying that $\nabla \cdot a = 0$.
For sufficiently small and smooth initial data, solutions exist globally in time and fulfill that
\begin{equation}\label{decayt}
	\| u(t) \|_{L^q (\mathbb{R}^n)} \le C t^{-\frac12} (1+t)^{-\gamma_q}
\end{equation}
for $1 \le q \le \infty$ and $\gamma_q = \frac{n}2 (1-\frac1q)$, and
\begin{equation}\label{decays}
	| u(t,x) | = O(|x|^{-n-1})
\end{equation}
as $|x| \to +\infty$ for any fixed $t$.
The details of these estimates can be found in \cite{AmrchGrltSchnbk,Brndls,Brndls-Vgnrn,Gg-Mykw-Osd,Kt,Mykw00}.
For other basic properties of Navier--Stokes flows, see \cite{FjtKt, KznOgwTnuch02,Lry} and references therein.
Particularly, \eqref{decayt} provides upper bound of the solution.
More detailed time global behavior of the solution is described by asymptotic approximations.
Asymptotic profiles have gained attention as a means of describing flow around structures such as obstructions, pumps and cylinders.
For this motivation, we refer to \cite{Glly-Mekw,Hshd,KznTrswWksg}.
Here we study quiet flow in the structure-free space.
For \eqref{NS} Fujigaki and Miyakawa \cite{Fjgk-Mykw} derived the asymptotic expansion of $u$ with $n$th order.
More precisely, they proved that there are unique smooth functions $U_m = U_m (t,x)$ such that $\lambda^{n+m} U_m (\lambda^2 t,\lambda x) = U_m (t,x)$ for $\lambda > 0$, and
\begin{equation}\label{asymplow}
	\biggl\| u(t) - \sum_{m=1}^n U_m (t) \biggr\|_{L^q (\mathbb{R}^n)}
	=
	O( t^{-\gamma_q - \frac{n}2 - \frac12} \log t)
\end{equation}
as $t \to +\infty$ for $1 \le q \le \infty$ under the condition that $(1+|x|)^{n+1} a \in L^1 (\mathbb{R}^n)$.
These $U_m$ are given as the concrete functions.
The logarithmic evolution in the estimate naturally emerges from the scale structure of the nonlinear advection.
We consider whether this logarithm is essential or not.
To determine this, a higher order expansion is required.
Such an expansion is derived from the Escobedo--Zuazua \cite{EZ} method together with the renormalization.
For the applications and the general theory of renormalization, see \cite{KtM,Ymd} and \cite{Ishg-Kwkm}, respectively.
On this way, the author \cite{Ymmt24,Ymmt25} identified the logarithmic evolution.

Escobedo--Zuazua method was developed to derive the asymptotic profiles of solutions to the convection-diffusion equation.
Here, we refer the bilinear equation that
\begin{equation}\label{cd}
\left\{
\begin{aligned}
	&\partial_t u - \Delta u = d\cdot \nabla (|u|u),
	&
	t>0,~ x \in \mathbb{R}^n,\\
	&u(0,x) = u_0 (x),
	&
	x \in \mathbb{R}^n,
\end{aligned}
\right.
\end{equation}
where $d \in \mathbb{R}^n$ is a constant.
They showed that solutions of \eqref{cd} have their own profiles $U_m$ such that $\lambda^{n+m} U_m (\lambda^2 t,\lambda x) = U_m (t,x)$ for $\lambda > 0$, and
\[
	\biggl\| u(t) - \sum_{m=1}^{n-2} U_m (t) \biggr\|_{L^q (\mathbb{R}^n)}
	=
	O( t^{-\gamma_q - \frac{n}2 + \frac12} \log t)  
\]
as $t \to+\infty$.
Especially, in two dimensions, this estimate is crucial and there exists $K_1\neq 0$ such that $\lambda^3 K_1 (\lambda^2 t, \lambda x) = K_1 (t,x)$ for $\lambda > 0$, and
\begin{equation}\label{asymp2dimcd}
    u(t) \sim U_0(t) + K_1(t) \log t
\end{equation}
as $t\to+\infty$.
Comparing this estimate and \eqref{asymplow} suggests that the nonlinearity of the Navier--Stokes flow is two steps weaker than one of \eqref{cd}.
This prediction is also supported by a comparison of scales.
Indeed, the condition $\nabla\cdot a = 0$ in \eqref{NS} is comparable to $\int_{\mathbb{R}^n} u_0 (x) dx = 0$ in \eqref{cd} and two equations are bilinear.
Furthermore, Carpio \cite{Crpo} proved the similar estimate as \eqref{asymp2dimcd} for two-dimensional Navier--Stokes flow with $\nabla\cdot a \neq 0$.
This implies that the effects of fluid outflow and inflow at the initial time develop logarithmically over time.

Let us return to the discussion of the Navier--Stokes equations.
To clarify whether \eqref{asymplow} is crucial or not, we derived the asymptotic expansion up to $2n$th order.
Namely, there exist unique smooth functions $U_m$ and $K_m$ such that $\lambda^{n+m} (U_m, K_m) (\lambda^2 t,\lambda x) = (U_m,K_m) (t,x)$ for $\lambda > 0$, and
\begin{equation}\label{simpre}
	u(t) \sim \sum_{m=1}^{2n} U_m (t) + \sum_{m=n+1}^{2n} K_m (t) \log t
\end{equation}
as $t \to +\infty$.
When the space dimension $n$ is even, this expansion is in optimal shape.
In particular, the logarithmic estimate \eqref{asymplow} is essential.
Since such logarithmic evolution is not found in the solution of linear problems, it can be read from \eqref{simpre} that nonlinear effects are strongly at work in even dimensional cases.
However, in the case $n$ is odd, logarithmic evolution never appear.
In fact,
\[
	u(t) \sim \sum_{m=1}^{2n} U_m (t)
\]
as $t \to +\infty$.
Hence the logarithmic estimate on \eqref{asymplow} is a pretense in the odd-dimensional cases.
Comparing these results, it can be seen that the characteristics of nonlinear effects appear weaker in the odd-dimension cases than in the even-dimension cases.
To state our assertion, we introduce the vorticity $\omega^{ij} = \partial_i u^j - \partial_j u^i$.
The vorticity satisfies that
\begin{equation}\label{v}
	\partial_t \omega^{ij} - \Delta \omega^{ij} +  \partial_i \mathcal{I}^j [u] - \partial_j \mathcal{I}^i [u] = 0,
\end{equation}
where
\[
	\mathcal{I}^j [u] = \omega^{\star j} \cdot u
\]
for $\omega^{\star j} = (\omega^{1j},\omega^{2j},\ldots,\omega^{nj})$.
We note that $\int_{\mathbb{R}^n} \mathcal{I}[u] (t,x) dx = 0$ (see \cite{Ymmt24,Ymmt25}).
To solve the Cauchy problem, we give the initial data $\omega (0) = \omega_0$ by $\omega_0^{ij} = \partial_i a^j - \partial_j a^i$.
As a natural conclusion, we expect that
\begin{equation}\label{lowfreqog}
    \int_{\mathbb{R}^n} x^\alpha \omega_0 (x) dx = 0
\end{equation}
for $|\alpha| \le 1$.
Indeed, we see that $\int_{\mathbb{R}^n} \omega_0^{ij} dx = \int_{\mathbb{R}^n} (\partial_i a^j - \partial_j a^i) dx = 0$ if $\omega_0$ and $a$ are in $L^1 (\mathbb{R}^n)$, and $\int_{\mathbb{R}^n} x_k \omega_0^{ij} dx = \int_{\mathbb{R}^n} (\delta_{jk} a^i - \delta_{ik} a^j) dx = 0$ from the integration by parts and the mass conservation law $\nabla \cdot a = 0$ (cf. \cite{Mykw00}), where $\delta_{jk}$ is Kronecker delta.
These conditions guarantee the decay estimates
\begin{equation}\label{decayog}
	\| \omega (t) \|_{L^q (\mathbb{R}^n)}
	\le
	C(1+t)^{-\gamma_q-1}
\end{equation}
for $1 \le q \le \infty$ and $\gamma_q = \frac{n}2 (1-\frac1q)$, and Biot--Savart law
\begin{equation}\label{BS}
	u^j = - \nabla (-\Delta)^{-1} \cdot \omega^{\star j},
\end{equation}
where $\omega^{\star j}$ is the $j$th column.
%
Moreover, when $\omega_0$ is localized, the estimates with weight
\begin{equation}\label{estwtog}
	\| |x|^k \omega (t) \|_{L^q (\mathbb{R}^n)}
	\le
	C t^{-\gamma_q} (1+t)^{-1+\frac{k}2}
\end{equation}
fulfill for some $k \in \mathbb{Z}_+$.
For details of these estimates, see \cite{Glly-Wyn06,Gg-Mykw-Osd,Kkvc-Trrs}.
Particularly, by Kukavica and Reis \cite{Kkvc-Ris}, these estimates are applied to consider the spatial decay of the velocity.
We study large-time behavior.
The assertion in the previous work is written as follows.
\begin{theorem}[cf. \cite{Ymmt24,Ymmt25}]\label{thm-pre}
Let the space dimension $n \ge 2$ be even, $\omega_0 \in L^1 (\mathbb{R}^n) \cap L^\infty (\mathbb{R}^n),~ |x|^{2n+1} \omega_0 \in L^1 (\mathbb{R}^n)$ and satisfy \eqref{lowfreqog} for $|\alpha| \le 1$.
Assume that the solutions $u$ of \eqref{NS} for $a^j = - \nabla (-\Delta)^{-1} \cdot \omega_0^{\star j}$ and $\omega$ of \eqref{v} for $\omega (0) = \omega_0$ meet \eqref{decayt}, \eqref{decayog} and \eqref{estwtog} for $k = 2n+1$, respectively, where $\omega_0^{\star j}$ is the $j$th column.
Then there exist unique functions $U_m$ and $K_m \in C((0,\infty),L^1 (\mathbb{R}^n) \cap L^\infty (\mathbb{R}^n))$ such that
\[
	\lambda^{n+m} (U_m,K_m) (\lambda^2 t,\lambda x) = (U_m,K_m) (t, x)
\]
for $\lambda > 0$, and
\begin{equation}\label{assertone}
	\bigg\| u(t) - \sum_{m=1}^{2n} U_m (t) - \sum_{m=n+1}^{2n} K_m (t) \log t \bigg\|_{L^q (\mathbb{R}^n)}
	= o(t^{-\gamma_q-n})
\end{equation}
as $t \to +\infty$ for $1 \le q \le \infty$ and $\gamma_q = \frac{n}2(1-\frac1q)$.
In addition, if $|x|^{2n+2} \omega_0 \in L^1 (\mathbb{R}^n)$, then the left-hand side of \eqref{assertone} is estimated by $O(t^{-\gamma_q-n-\frac12} (\log t)^2)$ as $t \to +\infty$.
\end{theorem}
These $U_m$ and $K_m$ will be provided as concrete functions.
Originally, this theorem was introduced for any dimensions, including odd dimensions other than one.
Examining the structure of $K_m$ again in detail, our main result is established as follows.
\begin{theorem}\label{thm-main}
Let the space dimension $n \ge 3$ be odd.
Assume the same conditions as in Theorem \ref{thm-pre} except for the space dimension.
Then there exist unique functions $U_m \in C((0,\infty),L^1 (\mathbb{R}^n) \cap L^\infty (\mathbb{R}^n))$ such that
\[
	\lambda^{n+m} U_m (\lambda^2 t,\lambda x) = U_m (t,x)
\]
for $\lambda > 0$, and
\begin{equation}\label{asserttwo}
	\biggl\| u(t) - \sum_{m=1}^{2n} U_m (t) \biggr\|_{L^q (\mathbb{R}^n)}
	= o(t^{-\gamma_q-n})
\end{equation}
as $t \to +\infty$ for $1 \le q \le \infty$ and $\gamma_q = \frac{n}2(1-\frac1q)$.
In addition, if $|x|^{2n+2} \omega_0 \in L^1 (\mathbb{R}^n)$, then the left-hand side of \eqref{asserttwo} is estimated by $O(t^{-\gamma_q-n-\frac12} \log t)$ as $t \to +\infty$.
\end{theorem}
In short, in odd-dimensional cases, the logarithmic evolutions $K_m$ disappear and the influence of the logarithm in the sharp estimate is mitigated.
Particularly, \eqref{asymplow} is not optimal.
Here, the condition $|x|^{n+2} \omega_0 \in L^1 (\mathbb{R}^n)$ is sufficient to show this fact.
This condition is compatible with $|x|^{n+1} a \in L^1 (\mathbb{R}^n)$.
Similarly $|x|^{2n+1}  \omega_0 \in L^1 (\mathbb{R}^n)$ corresponds to $|x|^{2n} a \in L^1 (\mathbb{R}^n)$.
The subtle differences in the second assertions should be noted, which the author expects to be optimal (see also Remark \ref{rmBrgrs}).
Like the previous theorem, this theorem requires no symmetry for the initial data.

\begin{remark}
After reading the proof, the reader might think that the conditions are given to the initial vortex for technical reasons.
Sure, when evaluating the velocity, conditions should be placed on the initial velocity.
In fact, this strategy is taken based on a more fundamental motivation.
Even if the initial velocity is localized, the velocity slowly decays as \eqref{decays}.
That is, $|x|^k u(t) \not\in L^1 (\mathbb{R}^n)$ for $k \in \mathbb{N}$ generally.
Is it reasonable to impose a property on the initial condition that the solution cannot satisfy?
This is also a problem related to global extensibility of time of Navier--Stokes flows.
On the other hand, $|x|^k \omega (t) \in L^1 (\mathbb{R}^n)$ could be guaranteed whenever the initial vorticity is localized.
\end{remark}

\begin{remark}
Why are logarithmic evolutions getting so much attention?
This serves as a marker for measuring the nonlinear effects on solutions.
From this perspective, \eqref{asymplow} can be considered an estimate that does not incorporate nonlinear characteristics.
For another reason why \eqref{asymplow} is a linear estimate, refer to the appendix.
\end{remark}

\begin{remark}\label{rmBrgrs}
One may predict that the second assertion of Theorem \ref{thm-main} is not optimal.
The equivalent problem to the one-dimensional Navier--Stokes equation is the Burgers equation with linear wave dissipation.
Namely, we refer to
\[
\left\{
\begin{aligned}
	&\partial_t u + u \partial_x u = \partial_x^2 u,
	&
	t>0,~ x \in \mathbb{R},\\
	&u(0,x) = a (x),
	&
	x \in \mathbb{R}\,
\end{aligned}
\right.
\]
for the initial-data such as $\int_{\mathbb{R}} a(x) dx = 0$.
If the second assertion of Theorem \ref{thm-main} is crutial, then it is expected that there exist $U_m$ and $K_3$ which fulfill $\lambda^{1+m} U_m (\lambda^2 t,\lambda x) = U_m (t, x)$ and $\lambda^4 K_3 (\lambda^2 t,\lambda x) = K_3 (t,x)$ and $u(t) \sim U_1 (t) + U_2 (t) + K_3 (t) \log t + U_3 (t)$ as $t \to +\infty$.
In fact, $K_3$ is vanishing.
\end{remark}

\paragraph{\textbf{Notations.}}
For a vector $u$ and a tensor $\omega$, we denote their components $j$th and $ij$th by $u^j$ and $\omega^{ij}$, respectively.
We abbreviate the $j$th column of $\omega$ by $\omega^{\star j} = (\omega^{1j},\omega^{2j},\ldots,\omega^{nj})$.
For vector fields $f$ and $g$, the convolution of them is simply denoted by $f*g (x) = \int_{\mathbb{R}^n} f(x-y) \cdot g (y) dy = \int_{\mathbb{R}^n} f(y) \cdot g(x-y) dy$.
We often omit the spatial parameter from functions, for example, $u(t) = u(t,x)$.
In particular, $G(t) * \omega_0 = \int_{\mathbb{R}^n} G(t,x-y) \omega_0 (y) dy$ and $\int_0^t g(t-s) * f(s) ds = \int_0^t \int_{\mathbb{R}^n} g(t-s,x-y) f(s,y) dyds$.
We symbolize the derivations by $\partial_t = \partial/\partial t,~ \partial_j = \partial/\partial x_j$ for $1 \le j \le n,~ \nabla = (\partial_1,\partial_2,\ldots,\partial_n)$ and $\Delta = \lvert \nabla \rvert^2 = \partial_1^2 + \partial_2^2 + \cdots + \partial_n^2$.
The length of a multiindex $\alpha = (\alpha_1,\alpha_2,\ldots, \alpha_n) \in \mathbb{Z}_+^n$ is given by $\lvert \alpha \rvert = \alpha_1 + \alpha_2 + \cdots + \alpha_n$, where $\mathbb{Z}_+ = \mathbb{N} \cup \{ 0 \}$.
We abbreviate that $\alpha ! = \alpha_1 ! \alpha_2! \cdots \alpha_n !,~ x^\alpha = x_1^{\alpha_1} x_2^{\alpha_2} \cdots x_n^{\alpha_n}$ and $\nabla^\alpha = \partial_1^{\alpha_1} \partial_2^{\alpha_2} \cdots \partial_n^{\alpha_n}$.
We define the Fourier transform and its inverse by $\hat{\varphi} (\xi) = \mathcal{F} [\varphi] (\xi) = (2\pi)^{-n/2} \int_{\mathbb{R}^n} \varphi (x) e^{-ix\cdot\xi} dx$ and $\check{\varphi} (x) = \mathcal{F}^{-1} [\varphi] (x) = (2\pi)^{-n/2} \int_{\mathbb{R}^n} \varphi (\xi) e^{ix\cdot\xi} d\xi$, respectively, where $i = \sqrt{-1}$.
The Riesz transforms are defined by $\mathcal{R}^j \varphi = \partial_j (-\Delta)^{-1/2} \varphi$, i.e., $\mathcal{R}^j \varphi =\mathcal{F}^{-1} [i\xi_j \lvert \xi \rvert^{-1} \hat{\varphi}]$ for $1 \le j \le n$ and $\mathcal{R} = (\mathcal{R}^1,\mathcal{R}^2,\ldots,\mathcal{R}^n)$.
Analogously, $\nabla (-\Delta)^{-1} \varphi = \mathcal{F}^{-1} [i\xi |\xi|^{-2} \hat{\varphi}]$.
The Lebesgue space and its norm are denoted by $L^q (\mathbb{R}^n)$ and $\| \cdot \|_{L^q (\mathbb{R}^n)}$, that is, $\| f \|_{L^q (\mathbb{R}^n)} = (\int_{\mathbb{R}^n} |f(x)|^q dx)^{1/q}$ for $1 \le q < \infty$ and $\| f \|_{L^\infty (\mathbb{R}^n)}$ is the essential supremum.
The heat kernel and its decay rate on $L^q (\mathbb{R}^n)$ are symbolized by $G(t,x) = (4\pi t)^{-n/2} e^{-|x|^2/(4t)}$ and $\gamma_q = \frac{n}2 (1-\frac1q)$.
We denote the floor function by Gauss symbol $[\mu]=\max\{ m \in \mathbb{Z} \mid m \le \mu \}$.
We employ Landau symbol.
Namely $f(t) = o(t^{-\mu})$ and $g(t) = O(t^{-\mu})$ mean $t^\mu f(t) \to 0$ and $t^\mu g(t) \to c$ for some $c \in \mathbb{R}$ such as $t \to +\infty$ or $t \to +0$.
When $\| u(t) - U(t) \|_{L^q (\mathbb{R}^n)} = o (\| U(t) \|_{L^q (\mathbb{R}^n)}$ for some $U$ and $1 \le q \le \infty$, we denote that $u(t)\sim U(t)$ as $t\to+\infty$.
We call an $n$-dimensional symmetric function $f$ `odd-type' when $f(-x)=-f(x)$.
Various positive constants are simply denoted by $C$.

\section{Preliminaries}
To prove our main result, we study the lower-order estimate \eqref{asymplow} in this section.
Using the vorticity, the mild solution of velocity is written as
\begin{equation}\label{MSu}
\begin{aligned}
	u^j (t)
	&=
	- \nabla (-\Delta)^{-1} G(t) * \omega_0^{\star j}
	-
	\int_0^t
		\mathcal{R}^j \mathcal{R} G (t-s) * \mathcal{I} [u] (s)
	ds
	-
	\int_0^t
		G(t-s) * \mathcal{I}^j [u] (s)
	ds,
\end{aligned}
\end{equation}
where 
$\mathcal{R} = (\mathcal{R}^1,\mathcal{R}^2,\ldots,\mathcal{R}^n)$ are Riesz transforms and $\mathcal{I} = (\mathcal{I}^1,\mathcal{I}^2,\ldots,\mathcal{I}^n)$ for $\mathcal{I}^j = \omega^{\star j} \cdot u$.
As already introduced, $\omega^{\star j}$ is the $j$th column of $\omega$.
This shape comes from coupling of the mild solution \eqref{v} and Biot--Savart law \eqref{BS}.
We emphasize that the integral kernels lead $\| \nabla (-\Delta)^{-1} G(t) * \varphi \|_{L^q (\mathbb{R}^n)} + \| \mathcal{R}^j \mathcal{R}G(t) * \phi \|_{L^q (\mathbb{R}^n)} \le C t^{-\gamma_q-\frac12}$ for $1 \le q \le \infty$ and some localized functions $\varphi$ and $\phi$ satisfying $\int_{\mathbb{R}^n} x^\alpha\varphi dx = 0$ for $|\alpha| \le 1$ and $\int_{\mathbb{R}^n} \phi dx = 0$.
Hence, time decay \eqref{decayt} is natural.
The same principle also provides the spatial decay \eqref{decays}.
For this principle, see \cite{Sbt-Smz,Stin,Zmr}.
We find the specific shape of $U_m$ on \eqref{simpre} from \eqref{MSu} by Escobedo--Zuazua method, then we see
\begin{equation}\label{Umlow}
\begin{aligned}
	&U_m^j (t)
	=
	- \sum_{|\alpha| = m+1} \frac{\nabla^\alpha \nabla (-\Delta)^{-1} G(t)}{\alpha!} \cdot \int_{\mathbb{R}^n}
		(-y)^\alpha \omega_0^{\star j} (y)
	dy\\
	&- \sum_{2l+|\beta| = m} \frac{\partial_t^l \nabla^\beta \mathcal{R}^j \mathcal{R} G(t)}{l!\beta!} \cdot \int_0^\infty \int_{\mathbb{R}^n}
		(-s)^l (-y)^\beta \mathcal{I}[u] (s,y)
	dyds\\
	&- \sum_{2l+|\beta| = m} \frac{\partial_t^l \nabla^\beta G(t)}{l!\beta!} \int_0^\infty \int_{\mathbb{R}^n}
		(-s)^l (-y)^\beta \mathcal{I}^j [u] (s,y)
	dyds
\end{aligned}
\end{equation}
for $1 \le m \le n$.
These functions, while seemingly different from those found in previous studies \cite{Crpo,Fjgk-Mykw}, are in fact the same (cf. \cite[Appendix A]{Ymmt25}).
When we expand the term of initial-data as
\[
\begin{aligned}
    &- \nabla (-\Delta)^{-1} G(t) * \omega_0^{\star j}
    =
    - \sum_{|\alpha| = 2}^N \frac{\nabla^\alpha \nabla (-\Delta)^{-1} G(t)}{\alpha!} \cdot \int_{\mathbb{R}^n}
		(-y)^\alpha \omega_0^{\star j} (y)
	dy
    + r_{N}^0 (t),
\end{aligned}
\]
Taylor theorem with \eqref{lowfreqog} and the estimates for Riesz transform say that the remaining part satisfies
\[
    \| r_{N}^0 (t) \|_{L^q (\mathbb{R}^n)} = o (t^{-\gamma_q-\frac{N-1}2})
\]
as $t \to +\infty$.
This is estimated by $O (t^{-\gamma_q-\frac{N}2})$ when $|x|^{N+1} \omega_0 \in L^1 (\mathbb{R}^n)$ is supposed.
We emphasize that whether $U_m$ is a odd or even-type function in $x$ is determined by $m$.
Here, we classify the parities of $n$-dimensional symmetric functions $f$ into `odd-type functions' and `even-type functions' according to whether $f(-x) = - f(x)$ or $f(-x) = f(x)$.
If $f$ is odd-type, then $\int_{\mathbb{R}^n} f(x) dx = 0$.
Some functions may be determined to be even-type even if they are odd in some $x_j$.
For example, we classify $x_j x_k e^{-|x|^2}$  for some $1 \le j,k \le n$ as even-type.
Since the parity is kept through the Fourier transform, $\mathcal{R}^j \mathcal{R}^k G$ is classified as even-type.
Under the assumption $1 \le m \le n$, $U_m$ is odd-type if $m$ is odd, and is even-type when $m$ is even.
Indeed, since $\nabla (-\Delta)^{-1} G$ is odd, the parity of $\nabla^\alpha \nabla (-\Delta)^{-1} G$ for $|\alpha|=m+1$ is decided by $m$.
Also, the even-oddness of $\partial_t^l \nabla^\beta \mathcal{R}^j \mathcal{R}^k G$ for $2l+|\beta|=m$ is determined by $m$.
From this perspective, the asymptotic profiles $U_m$ in \eqref{asymplow} appear analogous to that in the linear problem.
We also remark that, as $m$ becomes larger, such a structure breaks down and $U_m$ for $m \ge n+1$ consists of both partities (see \ref{appndJm}).
Since $\omega^{ij} = \partial_i u^j - \partial_j u_i$, the approximation of $\omega^{ij}$ is given by $\Omega_m^{ij} = \partial_i U_{m-1}^j - \partial_j U_{m-1}^i$ for $2 \le m \le n+1$.
Thus, the parity of $\Omega_m$ is opposite to that of $U_{m-1}$ and is also decided by $m$.
According to these symmetries, we conclude that
\begin{equation}\label{vanbs}
    \int_{\mathbb{R}^n} y^\beta (\Omega_{m_1}^{\star j} \cdot U_{m_2}) (t,x) dx = 0
\end{equation}
for $t>0$ if $|\beta|+m_1+m_2$ is odd.
This property plays a crucial role in next section.

\section{Proof of main theorem}
In this section, we prove Theorem \ref{thm-main}.
For the details of first half, refer to \cite{Ymmt24,Ymmt25}.
At first glance, it appears that the nonlinear term contains the logarithmic evolutions.
Indeed, on the same way as in \cite{Ymmt25}, we expand the second term of \eqref{MSu} as
\begin{equation}\label{GI}
\begin{aligned}
	&\int_0^t \mathcal{R}^j \mathcal{R}G (t-s) * \mathcal{I}[u] (s) ds\\
	&=
	\sum_{2l+|\beta| = 1}^{2n} \frac{\partial_t^l \nabla^\beta \mathcal{R}^j \mathcal{R}G(t)}{l!\beta!} \cdot \int_0^t \int_{\mathbb{R}^n}
		(-s)^l (-y)^\beta\\
    &\hspace{20mm} \biggl( \Bigl( \mathcal{I}[u] - \sum_{p=2}^{2l+|\beta|+1} \mathcal{I}_p \Bigr) (s,y) - \mathcal{I}_{2l+|\beta|+2} (1+s,y) \biggr)
	dyds\\
	&+
	\sum_{m=n+1}^{2n} \int_0^t \int_{\mathbb{R}^n}
		\biggl( \mathcal{R}^j \mathcal{R}G (t-s,x-y)\\
    &\hspace{20mm} - \sum_{2l+|\beta|=1}^m \frac{\partial_t^l \nabla^\beta \mathcal{R}^j \mathcal{R}G (t,x)}{l!\beta!} (-s)^l (-y)^\beta \biggr)
		\cdot \mathcal{I}_{m+2} (s,y)
	dyds\\
	&+
	\sum_{2l+|\beta| = n+1}^{2n} \frac{\partial_t^l \nabla^\beta \mathcal{R}^j \mathcal{R}G(t)}{l!\beta!} \cdot \int_0^t \int_{\mathbb{R}^n}
		(-s)^l (-y)^\beta \mathcal{I}_{2l+|\beta|+2} (1+s,y)
	dyds
	+ r_{2n+1}^1 (t)
\end{aligned}
 \end{equation}
for
\[
\begin{aligned}
	&r_{2n+1}^1 (t)
    =
	\int_0^t \int_{\mathbb{R}^n} \biggl(
		\mathcal{R}^j \mathcal{R}G (t-s,x-y) - \sum_{2l+|\beta| = 1}^{2n} \frac{\partial_t^l \nabla^\beta \mathcal{R}^j \mathcal{R}G (t,x)}{l!\beta!} (-s)^l (-y)^\beta
	\biggr)\\
	&\hspace{30mm}
	\cdot \Bigl( \mathcal{I}[u] - \sum_{p=n+3}^{2n+2} \mathcal{I}_p \Bigr) (s,y)
	dyds,
\end{aligned}
\]
where $\mathcal{I}_p$ for $n+3 \le p \le 2n+2$ 
are given by
\begin{equation}\label{defIp}
	\mathcal{I}_p^j = \sum_{m=1}^{p-n-2} \Omega_{p-n-m}^{\star j} \cdot U_m
\end{equation}
for $U_m$ defined by \eqref{Umlow}
and $\Omega_{m}^{ij} = \partial_i U_{m-1}^j - \partial_j U_{m-1}^i$.
It has the parabolic scale that
\begin{equation}\label{scIp}
    \lambda^{n+p} \mathcal{I}_p (\lambda^2 t,\lambda x) =  \mathcal{I}_p (t, x)
\end{equation}
for $\lambda > 0$, and fulfills that $\int_{\mathbb{R}^n} \mathcal{I}_p (t,x) dx = 0$.
Since $U_m$ and $\Omega_m$ are approximations of $u$ and $\omega$, respectively, $\mathcal{I}_p^j$ constitute an approximation of $\mathcal{I}^j [u] = \omega^{\star j} \cdot u$.
Hence, the coefficient on the first part of \eqref{GI} is represented as 
\begin{equation}\label{GI1}
\begin{aligned}
	&\frac1{l!\beta!}\int_0^t \int_{\mathbb{R}^n}
		(-s)^l (-y)^\beta \biggl( \Bigl( \mathcal{I}[u] - \sum_{p=2}^{2l+|\beta|+1} \mathcal{I}_p \Bigr) (s,y) - \mathcal{I}_{2l+|\beta|+2} (1+s,y) \biggr)
	dyds\\
	&=
	\mathcal{P}_{l\beta} (t^{-\frac12}) + \rho_{l\beta} (t)
\end{aligned}
\end{equation}
for some $(2n-2l-|\beta|)$th order polynomial $\mathcal{P}_{l\beta}$ and $\rho_{l\beta}$ satisfying $|\rho_{l\beta} (t)| = o (t^{-n+l+\frac{|\beta|}2})$ as $t\to +\infty$.
This claim was proved by employing \eqref{asymplow} and \eqref{estwtog}.
For the concrete shapes of $\mathcal{P}_{l\beta}$ and $\rho_{l\beta}$ and the details of this fact, see \cite[Claim 3.1]{Ymmt25}.
For $n+1 \le m \le 2n$, we put
\begin{equation}\label{GI2}
\begin{aligned}
	&J_m^j (t)
    = \int_0^t \int_{\mathbb{R}^n}
		\biggl( \mathcal{R}^j \mathcal{R} G (t-s,x-y)
        - \sum_{2l+|\beta|=1}^m \frac{\partial_t^l \nabla^\beta \mathcal{R}^j \mathcal{R} G (t,x)}{l!\beta!} (-s)^l (-y)^\beta \biggr)\\
		&\hspace{20mm}\cdot \mathcal{I}_{m+2} (s,y)
	dyds
\end{aligned}
\end{equation}
for the second part of \eqref{GI}.
Taylor theorem and Lebesgue convergence theorem guarantee that $J_m$ is well-defined in $C((0,\infty),L^1 (\mathbb{R}^n) \cap L^\infty (\mathbb{R}^n))$ and \eqref{scIp} provides the parabolic scale that
\[
     \lambda^{n+m} J_{m} (\lambda^2 t,\lambda x) = J_{m} (t,x)
\]
for $\lambda > 0$.
For the details of this fact, see also \cite[Claim 3.3]{Ymmt25}.
The symmetry of $J_m$ is determined by $m$.
However,
the even-oddness of $m$ and the even-oddness of $J_m$ present opposite correspondences.
The reader may doubt that $J_m$ contains both of odd and even-type functions.
This question is solved in \ref{appndJm}.
The story so far is exactly same as in \cite{Ymmt24,Ymmt25}.
In those works, the third part of the right-hand side of \eqref{GI} were thought to provide the logarithmic evolutions because the spatiotemporal integrations are separated to
\[
\begin{aligned}
	&\int_0^t \int_{\mathbb{R}^n}
		(-s)^l (-y)^\beta \mathcal{I}_{2l+|\beta|+2} (1+s,y)
	dyds
	=
	\int_0^t s^l (1+s)^{-l-1} ds
	\int_{\mathbb{R}^n}
		(-1)^l (-y)^\beta \mathcal{I}_{2l+|\beta|+2} (1,y)
	dy
\end{aligned}
\]
for $n+1 \le 2l+|\beta| \le 2n$.
Here $\int_0^t s^l (1+s)^{-l-1} ds$ surely provides $\log t$.
Namely, from the third part of \eqref{GI}, we see
\begin{equation}\label{GI3}
\begin{aligned}
    &\sum_{2l+|\beta|=m} \frac{\partial_t^l \nabla^\beta \mathcal{R}^j \mathcal{R}G(t)}{l!\beta!} \cdot \int_0^t \int_{\mathbb{R}^n}
		(-s)^l (-y)^\beta \mathcal{I}_{2l+|\beta|+2} (1+s,y)
	dyds\\
    &=
    \tilde{K}_m^j (t) \log t
    + \sum_{2l+|\beta|=m} \mathcal{Q}_{l\beta} (t^{-1})\partial_t^l \nabla^\beta \mathcal{R}^j \mathcal{R}G(t)
    + \sum_{2l+|\beta|=m} \varrho_{l\beta} (t)\partial_t^l \nabla^\beta \mathcal{R}^j \mathcal{R}G(t)
\end{aligned}
\end{equation}
for
\[
    \tilde{K}_m^j (t) = \sum_{2l+|\beta|=m} \frac{\partial_t^l \nabla^\beta \mathcal{R}^j \mathcal{R}G(t)}{l!\beta!} \cdot \int_{\mathbb{R}^n}
		(-1)^l (-y)^\beta \mathcal{I}_{2l+|\beta|+2} (1,y)
	dy
\]
and some $[n-l-\frac{|\beta|}2]$th order polynomial $\mathcal{Q}_{l\beta}$ and $\varrho_{l\beta}$ satisfying $|\varrho_{l\beta} (t)|=o(t^{-n+l+\frac{|\beta|}2})$ as $t\to+\infty$.
The concrete shapes of $\mathcal{Q}_{l\beta}$ and $\varrho_{l\beta}$ are found in \cite[Claim 3.2]{Ymmt25}.
Substituting \eqref{GI1}-\eqref{GI3} to \eqref{GI} yields that
\begin{equation}\label{expNLfin}
\begin{aligned}
    &\int_0^t \mathcal{R}^j \mathcal{R}G (t-s) * \mathcal{I}[u] (s) ds\\
    &=
   \sum_{2l+|\beta|=1}^{2n} \left( \mathcal{P}_{l\beta}(t^{-1/2}) + \mathcal{Q}_{l\beta}(t^{-1}) \right) \partial_t^l\nabla^\beta \mathcal{R}^j\mathcal{R} G(t)
   +\sum_{m=n+1}^{2n} J_m^j (t)
   + \sum_{m=n+1}^{2n} \tilde{K}_m^j (t)\log t\\
   & + r_{2n+1}^1(t) + r_{2n+1}^2(t)
\end{aligned}
\end{equation}
for the above $r_{2n+1}^1$
and
\[
    r_{2n+1}^2 (t)
    =
    \sum_{2l+|\beta|=1}^{2n} \left( \rho_{l\beta}(t)+\varrho_{l\beta} (t) \right) \partial_t^l\nabla^\beta \mathcal{R}^j\mathcal{R} G(t).
\]
The decays of $\rho_{l\beta}$ and $\varrho_{l\beta}$ yield that
$
    \| r_{2n+1}^2 (t) \|_{L^q (\mathbb{R}^n)}
    = O(t^{-\gamma_q-n-\frac12})
$
as $t\to+\infty$.
Another term fulfills that
$
    \| r_{2n+1}^1 (t) \|_{L^q (\mathbb{R}^n)}
    = O(t^{-\gamma_q-n-\frac12}\log t)
$
as $t\to+\infty$.
These facts are also confirmed in \cite{Ymmt25}.
The third term of \eqref{MSu} is handled on the same way.
Rearranging these terms according to their parabolic scales, we finally decide $U_m$ satisfying $\lambda^{n+m} U_m (\lambda^2 t, \lambda x) = U_m (t,x)$ for $\lambda > 0$ and have that
\[
    u(t) = \sum_{m=1}^{2n} U_m (t) + \sum_{m=n+1}^{2n} K_m (t) \log t + r_{2n+1} (t)
\]
for
\[
\begin{aligned}
    &K_m^j (t) = - \sum_{2l+|\beta|=m} \frac{\partial_t^l \nabla^\beta \mathcal{R}^j \mathcal{R}G(t)}{l!\beta!} \cdot \int_{\mathbb{R}^n}
		(-1)^l (-y)^\beta \mathcal{I}_{2l+|\beta|+2} (1,y)
	dy\\
    &- \sum_{2l+|\beta|=m} \frac{\partial_t^l \nabla^\beta G(t)}{l!\beta!} \int_{\mathbb{R}^n}
		(-1)^l (-y)^\beta \mathcal{I}_{2l+|\beta|+2}^j (1,y)
	dy
\end{aligned}
\]
and $r_{2n+1}$ satisfying $\| r_{2n+1} (t) \|_{L^q (\mathbb{R}^n)} = O (t^{-\gamma_q-n-\frac12} \log t)$ as $t\to +\infty$.
In fact, from \eqref{vanbs} and \eqref{defIp}, we see for the coefficient of $K_m$ that
\begin{equation}\label{oddeven}
\begin{aligned}
	&\int_{\mathbb{R}^n}
		(-1)^l (-y)^\beta \mathcal{I}_{2l+|\beta|+2}^j (1,y)
	dy
    =
	\sum_{m=1}^{2l+|\beta|-n} \int_{\mathbb{R}^n} (-1)^l (-y)^\beta (\Omega_{2l+|\beta|+2-n-m}^{\star j} \cdot U_m)(1,y) dy = 0
\end{aligned}
\end{equation}
and then $K_m (t) = 0$.
Therefore, $u$ contains no logarithms and we complete the proof.

\begin{remark}
In even dimensions, it is difficult to determine whether $K_m$ diverges or not.
Futhermore, since $K_m$ consists of several terms, we should confirm whether they cancel each other or not.
\end{remark}

\appendix
\section{Asymmetry due to nonlinear effects}\label{appndJm}
We consider the parity of $J_m$.
We already saw that the even-oddness of $\mathcal{I}_p$ is determined by $p$.
However, for odd dimensions, $\mathcal{I}_p (s,-y) = - \mathcal{I}_p (s,y)$ when $p$ is even, and $\mathcal{I}_p (s,-y) = \mathcal{I}_p (s,y)$ if $p$ is odd.
Incidentally, for even dimensions, these correspondences are reversed.
Now we study odd dimensions.
Hence, the even-oddness of the first term $\mathcal{R}^j \mathcal{R} G(t-s) * \mathcal{I}_{m+2} (s)$ is decided.
Among the other integrals $\int_{\mathbb{R}^n} (-y)^\beta \mathcal{I}_{m+2} dy$, only the cases $|\beta|+m$ being odd can be remained  since any even-types are vanishing.
Therefore, the parity of $J_m$ and one of $m$ correspond in reverse order.
This is a distinctly different property from that which appears in the linear problem.
In other words, $J_m$ are the embodiments of the nonlinear effect.
Also, the term on the first part of \eqref{expNLfin} contains several parities.
The even-oddness in $x$ and the spatial decay-rate of this term are decided when $(l,\beta)$ is fixed.
On the other hand, this term contains several decay-rates in time.
Such an expansion never be seen in linear problems.
These terms reveal the distortion of the solution's symmetry due to nonlinear effects.




\end{document}